\documentclass{article}
\usepackage{xfrac}
\usepackage{textcomp}
\usepackage{amsmath,graphicx}

\begin{document}

\title{Visualization of Fractional Integrals}
\author{T.L. Grobler: Stellenbosch University - Computer Science Division.}

\maketitle

\begin{abstract}
We presented a novel geometric interpretation of the Riemann-Liouville fractional integral. We found that a Riemann-Liouville integral can 
be thought of as the area obtained by summing together the area of an infinite number of non-rectangular infinitesimals whose shape is determined by the 
order of integration $\alpha$ and the integration limit $t$. We also showed that this geometric interpretation offers many pedagogical benefits as it is very similar
in nature to the geometric interpretation of the Riemann integral.   
\end{abstract}



\noindent
Leibniz, who is regarded by many scholars as the father of calculus, is credited with inventing the following notation: $\frac{d^n y}{d x^n}$. 
A natural question which arises when one inspects this notation is: Can $n$ be a non-integer value? L'Hospital phrased it differently: ``\emph{What if $n$ be a $\frac{1}{2}$?}''. 
In a letter Leibniz sent to L'Hospital, Leibniz comments on L'Hospital's question by making the following prophetic remark: ``\emph{This is an apparent paradox from which, one day, useful 
consequences will be drawn}'' \cite{ross77,machado14,machado17}.\\

\noindent
The mathematical discipline in which we study the extension of derivatives and integrals to non-integer orders is known as fractional calculus. Euler, Liouville and Riemann laid most 
of the theoretical groundwork, which was needed to develop the field of fractional calculus \cite{davis59,euler1738,liouville1832,riemann1876,laurent1884}. Fractional calculus has now developed into a mature discipline \cite{machado14,machado17}. 
However, there is one aspect of the discipline which only started to develop fairly recently: its physical and geometric interpretability. At the first international conference on fractional calculus in
New Haven (USA) (which took place in 1974) it was pointed out that literature lacked an acceptable geometric and physical interpretation of fractional calculus \cite{ross06}. Subsequently, the last few decades have seen 
a number of papers that attempt to rectify this. There now exist probabilistic interpretations \cite{stanislavsky04,machado03}, geometric interpretations \cite{adda97,tarasov16,podlubny02}, physical interpretations \cite{cioc16,rutman95, nigmatullin92, gomez14} and economic 
interpretations \cite{tarasova17} of fractional order integrals and differentials.\\

\noindent
In this article we specifically focus on the geometric interpretation of fractional integrals. To be more specific, we present a novel geometric interpretation of the Riemann-Liouville fractional integral based at 0 \cite{laurent1884}.
As we mentioned earlier, Riemann and Liouville played an important role in the development of fractional calculus. Their work was crucial in defining what is now known as 
the Riemann-Liouville fractional integral \cite{laurent1884}. At this point in time the reader should be made aware of the fact that there exists many other accepted fractional integral definitions in the literature.
We, however, deem these definitions to be beyond the scope of the current article and as such we do not discuss them any further in this work. Furthermore, the geometric interpretation 
we present here is an extension of the geometric interpretation presented by Podlubny \cite{podlubny02}. Podlubny showed that a Riemann-Liouville fractional integral 
can be converted into a Riemann-Stieltjes integral. Furthermore, Grobler showed that a Riemann-Stieltjes integral can be interpreted as the area which is obtained by summing together 
the area of an infinite number of non-rectangular infinitesimals, i.e. a Riemann-Stieltjes integral can be converted into a Cavalieri integral \cite{ackermann12,grobler19}. Combining these two ideas results in the following geometric 
interpretation: the Riemann-Liouville fractional integral can be interpreted as the area which is obtained by summing together the area of an infinite number of non-rectangular infinitesimals whose shape is determined by the order of integration $\alpha$ and the integration limit $t$.\\

\noindent
It is important for the reader to realize that we specifically wrote this paper for people with a pedagogical interest, undergraduate students and for the layman who finds integration theory intriguing and 
as such we will try to steer clear from using too many mathematical definitions and proofs. We will rather convey our ideas using as many examples as possible. From the 
three groups mentioned we believe that the undergraduate student will find this article most useful. Although undergraduate students tend to have a limited understanding of 
fractional calculus they generally are quite interested in the subject. This is evident from perusing mathematical social media platforms and from surveying 
articles published in undergraduate journals \cite{munkhammar05}. Undergraduates struggle to grasp fractional calculus, because it is very abstract. We believe that the geometric interpretation we present in this paper, as it is similar in nature to the definition of classical integration, will provide undergraduates with the insight that they need to quickly grasp and understand the fundamental 
concepts that underpin fractional calculus.\\


\noindent
We start the paper by reviewing the definition of the Riemann-Stieltjes integral, the Cavalieri integral and the Euler function. We then present the definition of the Riemann-Liouville fractional integral and we end the paper by presenting our novel geometric interpretation, some examples and a conclusion.

\section{Riemann-Stieltjes Integral}
Recall that the Riemann-Stieltjes integral is defined as follows \cite{bartle76}:
\begin{equation}
\label{eq:g_int2}
\int_{a'}^{b'} f(x) dg(x) =  \lim_{n \rightarrow \infty}\sum_{i=0}^{n-1} f(x_i^2)[g(x_{i+1}^2)-g(x_{i}^2)], 
\end{equation}
\noindent
with $(x_i^2)_{i=0}^n$ being arbitrary partition points on the $x$-axis. 
The function $g$ is known as the integrator and if it is monotone increasing then it maps points in the interval $[a',b,]$ to the interval $[a,b]$, with $a = g(a')$ and $b = g(b')$.
The integral in Eq.~\eqref{eq:g_int2} is hard to evaluate in its present form. If $g$ is differentiable, however, then it becomes trivial to 
evaluate a Riemann-Stieltjes integral due to the following identity:
\begin{equation}
\int_{a'}^{b'} f(x) dg(x) = \int_{a'}^{b'} f(x)g'(x)~dx.
\end{equation}
As an example:
\begin{equation}
\int_{\frac{1}{2}}^{2} f(x) dg(x)=\int_{\frac{1}{2}}^2 x~d2x = \int_{\frac{1}{2}}^2 2x~dx = 3.75. 
\end{equation}

\section{Cavalieri Integral}
\begin{figure}[htb]
\centering
\includegraphics[width=0.75\textwidth]{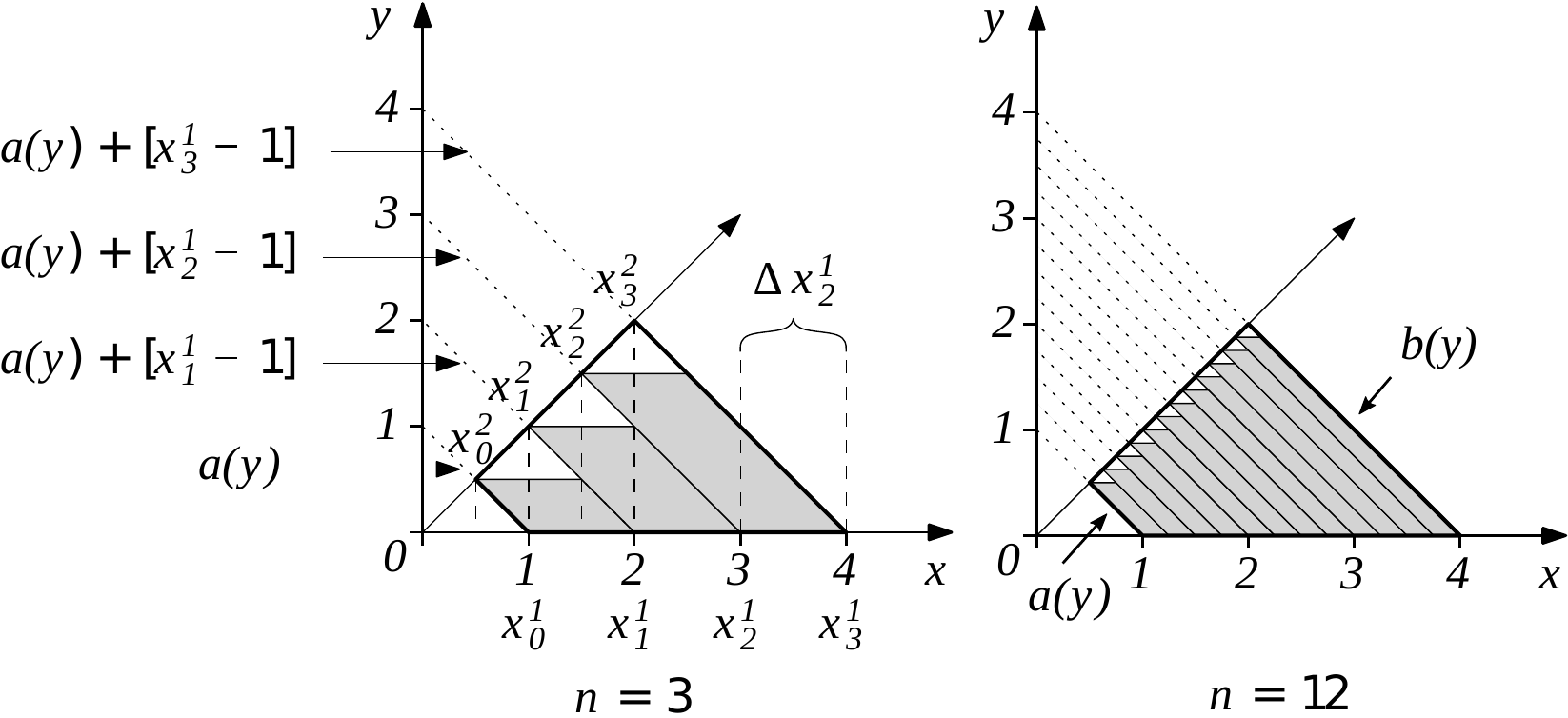}
\caption{Region bounded by the $x$-axis and the lines $f(x)=x$, $a(y)=1-y$, and $b(y)=4-y$. The figure also depicts the partition points $x_i^2$ as used in the Cavalieri sum (see Eq.~\eqref{eq:cav_sum}). Reproduced from Quaestiones Mathematicae (2012) 35: 265-296 with permission \copyright~ NISC (Pty) Ltd.}
\label{fig:caval2}
\end{figure}

In this section we define the Cavalieri integral \cite{ackermann12}. Let the region $R$ be bounded by the $x$-axis and the lines $f(x)=x$, $a(y)=1-y$ and $b(y)=4-y$. This region is shown in Fig.~\ref{fig:caval2}.\\

\noindent
The area of $R$ can be determined as follows if we employ the classical notion of integration:
\begin{equation}
\int_0^2x\, dx+\int_2^44-x\, dx- \int_0^{\frac{1}{2}}x\, dx-\int_{\frac{1}{2}}^11-x\, dx = 3.75. 
\end{equation}

\noindent
There, however, exists a more straightforward way do determine the area of $R$, i.e. we can sum together the area of non-rectangular integration strips inscribing $R$ instead of
having to sum together the area of multiple rectangular integration strips. Moreover, the shape of the non-rectangular integration strips we should use is determined by the function $a(y)$. We can express this notion more formally as 
follows, let $(x_i^1)_{i=0}^{n}$ denote a partition on the $x$-axis, such that $a = x_0^1 < x_1^1 < \cdots < x_n^1 = b$, and $\Delta x_i^1 = x_{i+1}^1 - x_i^1$.
We are now able to construct the following sum (the lower Cavalieri sum):
\begin{equation}
\label{eq:cav_sum}
\sum_{i=0}^{n-1} f(x_i^2)\Delta x_i^1.
\end{equation}
The partition points $(x_i^2)_{i=0}^{n}$ are depicted in Fig.~\ref{fig:caval2}. Note that Eq.~\eqref{eq:cav_sum} approximates the area of $R$. In the limit Eq.~\eqref{eq:cav_sum} approaches 
the Cavalieri integral:
\begin{equation}
\label{eq:caval1}
\int_{a(y)}^{b(y)}f(x)\, dx = \lim_{n\to \infty}\sum_{i=0}^{n-1} f(x_i^2)\Delta x_i^1.
\end{equation}
It is, however, quite hard to evaluate the above integral directly. Fortunately, it is easy to convert a Cavalieri integral into an equivalent Riemann or Riemann-Stieltjes 
integral by using the transformation functions $h$ and $g$. Expressed mathematically:
\begin{equation}
\label{eq:main_cav}
\int_{a(y)}^{b(y)}f(x)\,dx =\int_a^b f \circ h (x)\, dx = \int_{a'}^{b'} f(x) dg(x).
\end{equation}
We can calculate $g$ as follows:
\begin{equation}
g(x) = x - a\circ f(x) + a.
\end{equation}
Moreover, $h=g^{-1}$. We can now evaluate the Cavalieri integral by using its equivalent Riemann or Rieman-Stieltjes integral.
If we use its Riemann equivalent we obtain:
\begin{equation}
\int_{a(y)}^{b(y)}f(x)\, dx = \int_a^b f \circ h (x)\, dx = \dfrac{1}{2}\int_1^4x\, dx = 3.75.  
\end{equation}
If we use its Riemann-Stieltjes equivalent we obtain:
\begin{equation}
\int_{a(y)}^{b(y)}f(x)\, dx = \int_{a'}^{b'} f \, dg(x) = \int_{\frac{1}{2}}^2x\, d2x = 3.75.  
\end{equation}

\noindent
Conversely, if $g$ is known (and not $a(y)$) and $g(a') = a$ then we can calculate $a(y)$ using the following \cite{grobler19}:
\begin{equation}
\label{eq:a_y}
a(y) = f^{-1}(y) - g\circ f^{-1}(y) + g(a'). 
\end{equation}
Eq.~\eqref{eq:a_y} allows us to transform a Riemann-Stieltjes integral into an equivalent Cavalieri integral. 
Furthermore, this transformation enables us to assign a geometric interpretation to a Riemann-Stieltjes integral. The interpretation being:
it represents the area obtained by summing together the area of an infinite number of non-rectangular infinitesimals; all of them being translations of $a(y)$. Note,
$b(y) = a(y) + (b-a)$.

\section{Gamma Function}
\begin{figure}[htb]
\centering
\includegraphics[width=0.75\textwidth]{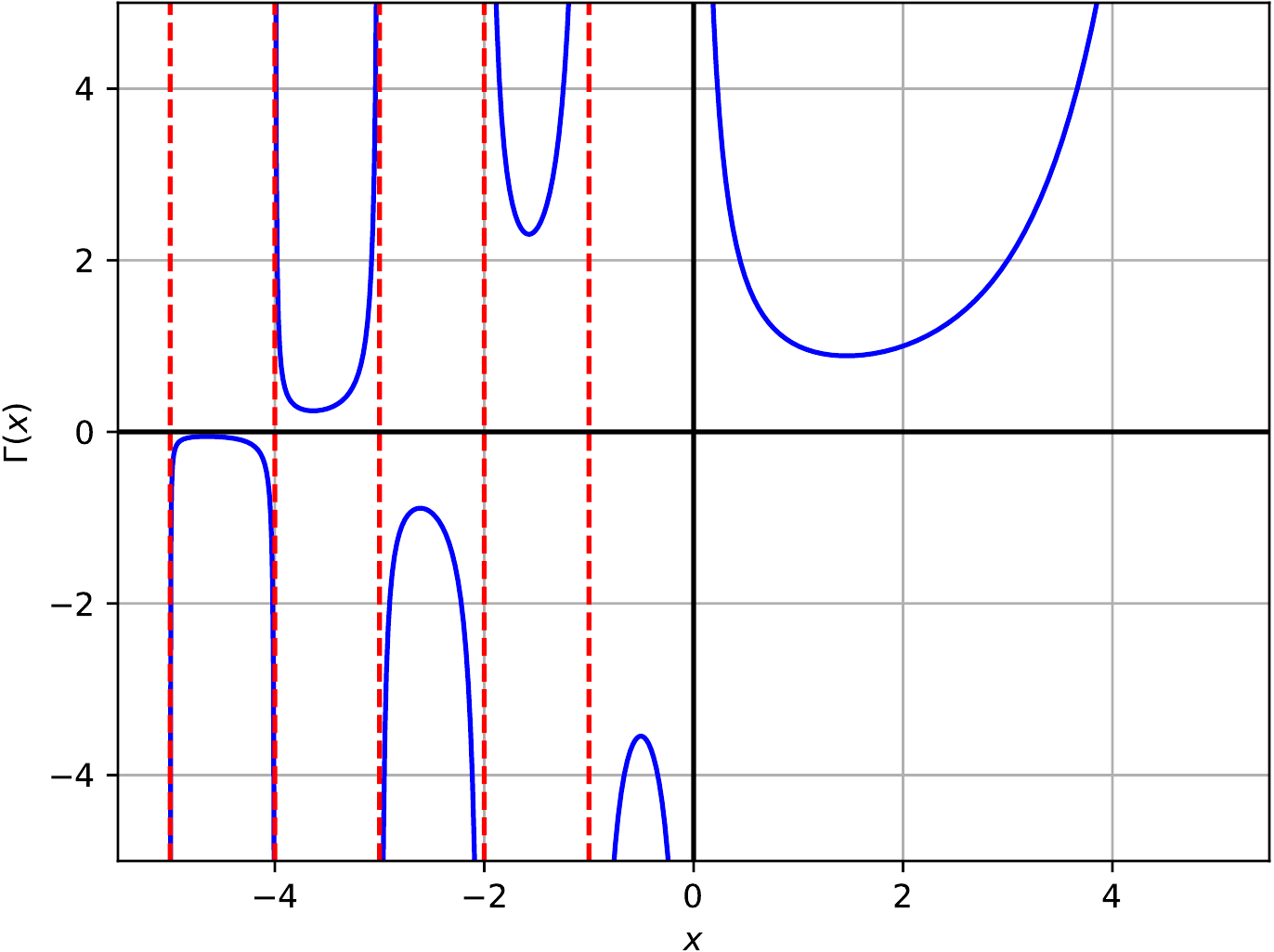}
\caption{A plot of $\Gamma(x)$ for $x\in(-5,5)$. $\Gamma(x)$ is not defined for the negative integers.}
\label{fig:gamma}
\end{figure}

The gamma function was first devised by Euler \cite{euler1738}. The gamma function is defined for all real numbers except for the negative integers. Moreover, it is a natural extension of the factorial function to the real numbers. The Gamma function is defined as:
\begin{equation}
\Gamma(x) = \int_0^{\infty} t^{x-1} e^{-t}\,dt.
\end{equation}
The function $\Gamma(x)$ is depicted in Fig.~\ref{fig:gamma} for $x\in(-5,5)$.

\section{Riemann-Liouville Fractional Integral}
In this section we present the definition of the Riemann-Liouville fractional integral of order $\alpha$ \cite{laurent1884}. Let
\begin{equation}
If(t) := \int_0^t f(\tau)~d\tau.
\end{equation}
If we apply the above operator in a repetitive manner to $f$ we obtain the $n$th antiderivative of $f$ based at 0:
\begin{equation}
\label{eq:n_anti}
I^nf(t) = \int_0^t\int_0^{t_1}\cdots \int_0^{t_{n-1}}f(t_n)~dt_n\cdots dt_2 dt_1.
\end{equation}
Applying Cauchy's formula for repetitive integration to the above $n$th antiderivative of $f$ helps us express the above $n$th antiderivative as a single integral. 
Applying Cauchy's formula to Eq.~\eqref{eq:n_anti} results in:
\begin{equation}
\label{eq:cauchy}
I^nf(t) = \frac{1}{(n-1)!}\int_0^t (t-\tau)^{n-1}f(\tau)~d\tau.
\end{equation}
As an example if $f(\tau)=\tau$ and $n=2$ then Eq.~\eqref{eq:cauchy} reduces to:
\begin{equation}
I^2f(t) = \int_0^t (t-\tau)t~dt = \frac{\tau^2 t}{2} - \frac{\tau^3}{3} \Bigg |_0^t = \frac{t^3}{3!}.
\end{equation}
The definition in Eq.~\eqref{eq:cauchy} can now be extended to an arbitrary fractional order by replacing $(n-1)!$ with the Gamma function (see Fig.~\ref{fig:gamma}):
\begin{equation}
\label{eq:frac_integral}
I^{\alpha}f(t) = \frac{1}{\Gamma(\alpha)}\int_0^t (t-\tau)^{\alpha-1}f(\tau)~d\tau.
\end{equation}
Eq.~\eqref{eq:frac_integral} is known as the Riemann-Liouville fractional integral of order $\alpha$. As an example if $f(\tau)=\tau$ and $\alpha=\frac{1}{2}$ then 
Eq.~\eqref{eq:frac_integral} reduces to:
\begin{equation}
I^{\frac{1}{2}}f(t) = \frac{1}{\Gamma(\frac{1}{2})} \int_0^t \frac{\tau}{\sqrt{t-\tau}} dt = \frac{1}{\sqrt{\pi}}\left [ -\frac{2}{3}\sqrt{t-\tau}(t+2x)\right] \Bigg |_0^t=\frac{4}{3\sqrt{\pi}}t^{\frac{3}{2}}. 
\end{equation}
Furthermore, it is straightforward to show that the $I$ operator satisfies:
\begin{equation}
I^{\alpha}I^{\beta}f(t) = I^{\alpha+\beta}f(t). 
\end{equation}

\section{Geometric Interpretation of the Riemann-Liouville Fractional Integral}
In this section we present a novel geometric interpretation of the Riemann-Liouville fractional integral. 
Recall that the Riemann-Liouville fractional integral of order $\alpha$ is defined as:
\begin{equation}
\label{eq:RL_frac_int}
\frac{1}{\Gamma(\alpha)}\int_{0}^{t}f(\tau)(t-\tau)^{\alpha-1}~d\tau
\end{equation}

\noindent
Podlubny showed that the above fractional integral can be rewritten as a Riemann-Stieltjes integral \cite{podlubny02}:
\begin{equation}
\label{eq:rs}
\int_0^{t} f(\tau)~dg_t(\tau),
\end{equation}
with 
\begin{equation}
\label{eq:g_rl}
g_t^{\alpha}(\tau) = \frac{\left \{t^{\alpha} - (t-\tau)^{\alpha} \right \}}{\Gamma(\alpha+1)}. 
\end{equation}

\noindent
Note that $g_t^{\alpha}(\tau)$ is always invertible, since $\frac{d g_t^{\alpha}}{d \tau} = \frac{(t-\tau)^{\alpha-1}}{\Gamma(\alpha)}>0$ for $\tau\in[0,t]$. The inverse 
of $g_t^{\alpha}(\tau)$ is equal to:
\begin{equation}
h_t^{\alpha}(\tau) = t - \sqrt[\alpha]{t^{\alpha} - \Gamma(\alpha+1)\tau}
\end{equation}

\begin{figure}[htb]
\centering
\includegraphics[width=0.75\textwidth]{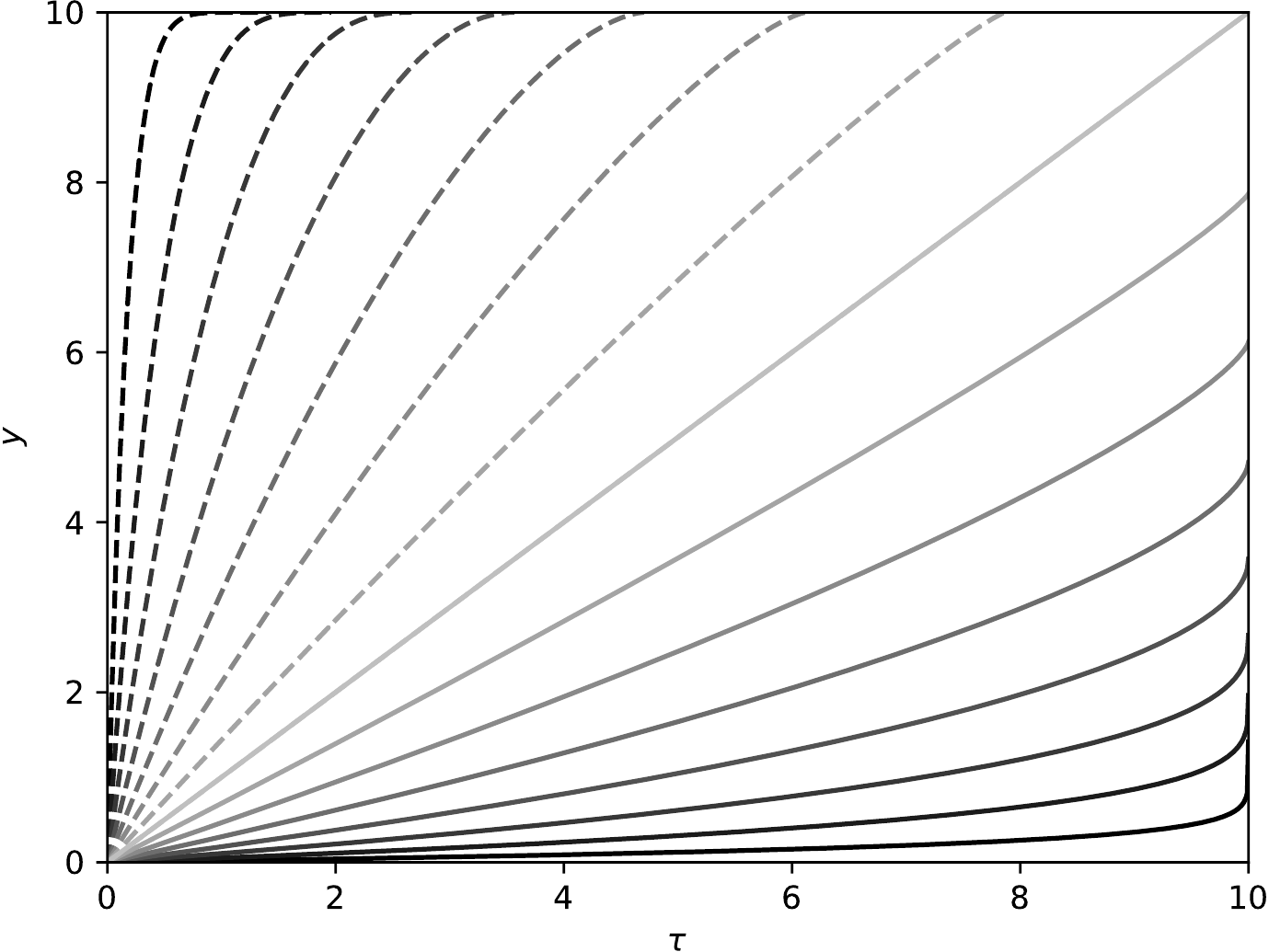}
\caption{This figure depicts the curves of $g_{10}^{\alpha}(\tau)$ and $h_{10}^{\alpha}(\tau)$ for $\alpha\in\mathsf{A}$. The curves associated with $g$ are plotted using 
solid lines, while the curves plotted with dashed lines are associated with $h$. The lighter the shade with which a curve is depicted, the larger the $\alpha$-value is that is associated with it.
}
\label{fig:gandh}
\end{figure}

\noindent
Let $\mathsf{A} = \{0,0.2,0.4,0.6,0.8,1\}$. The analytic expressions of $g_{10}^{\alpha}(\tau)$ and $h_{10}^{\alpha}(\tau)$ are presented in Table~\ref{tab:gandh} for $\alpha\in\mathsf{A}$. 
The expressions in Table~\ref{tab:gandh} are depicted in Fig~\ref{fig:gandh}.

\begin{table}[h!]
 \centering
 \caption{The analytic expressions of $g_{10}^{\alpha}(\tau)$ and $h_{10}^{\alpha}(\tau)$ for $\alpha\in\mathsf{A}$.}
 \label{tab:gandh}
 \begin{tabular}{|c || c | c|} 
 \hline
 $\alpha$ & $g_{10}^{\alpha}(\tau)$ & $h_{10}^{\alpha}(\tau)$ \\ [0.5ex] 
 \hline\hline
 $0$ & 0 & $y=0$ \\ 
 $0.2$ & $[\Gamma\left (\frac{6}{5} \right )]^{-1}\left(\sqrt[5]{10}-\sqrt[5]{10-\tau}\right)$ & $10 - \left ( \sqrt[5]{10} -  \Gamma\left (\frac{6}{5} \right ) \tau \right )^5$  \\
 $0.4$ & $[\Gamma\left (\frac{7}{5} \right )]^{-1}\left(\sqrt[5]{10^2}-\sqrt[5]{(10-\tau)^2}\right)$ & $10 - \sqrt{\left ( \sqrt[5]{10} -  \Gamma\left (\frac{7}{5} \right ) \tau \right )^5}$ \\
 $0.6$ & $[\Gamma\left (\frac{8}{5} \right )]^{-1}\left(\sqrt[5]{10^3}-\sqrt[5]{(10-\tau)^3}\right)$ & $10 - \sqrt[3]{\left ( \sqrt[5]{10} -  \Gamma\left (\frac{8}{5} \right ) \tau \right )^5}$ \\
 $0.8$ & $[\Gamma\left (\frac{9}{5} \right )]^{-1}\left(\sqrt[5]{10^4}-\sqrt[5]{(10-\tau)^4}\right)$ & $10 - \sqrt[4]{\left ( \sqrt[5]{10} -  \Gamma\left (\frac{9}{5} \right ) \tau \right )^5}$ \\ [1ex] 
 $1$ & $\tau$ & $\tau$ \\ [1ex] 
 \hline
 \end{tabular}
 \end{table}

\noindent
Furthermore, Grobler showed that it is trivial to convert a Riemann-Stieltjes integral into a Cavalieri integral \cite{ackermann12,grobler19} (see Fig.~\ref{fig:caval2} and the Cavalieri section of this paper if more details are required). If we apply Grobler's conversion method to Eq.~\eqref{eq:rs} we obtain:
\begin{equation}
\label{eq:cav}
\int_{a_t^{\alpha}(y)}^{b_t^{\alpha}(y)} f(\tau)~d\tau, 
\end{equation}
with
\begin{equation}
\label{eq:cav_a}
a_t^{\alpha}(y) = f^{-1}(y) - \frac{t^{\alpha}-(t-f^{-1}(y))^{\alpha}}{\Gamma(\alpha+1)}.
\end{equation}
Eq.~\eqref{eq:cav_a} is obtained by substituting Eq.~\eqref{eq:g_rl} into Eq.~\eqref{eq:a_y}. Note, $b_t^{\alpha}(y) = a_t^{\alpha}(y) + \frac{t^{\alpha}}{\Gamma(\alpha+1)}$. 

\noindent
We can now assign a geometric interpretation to Eq.~\eqref{eq:cav}, since it is a Cavalieri integral (see Fig.~\ref{fig:caval2} and the Cavalieri section of this paper). The fractional integral in Eq.~\eqref{eq:RL_frac_int} can, therefore, be interpreted as the area obtained 
by summing together the area of an infinite number of infinitesimally small non-rectangular integration strips whose shape is determined by $\alpha$ and $t$. If $\alpha$ is equal to one, however, then the integral reduces to a normal Riemann integral as the integration strips become rectangular for this particular choice of $\alpha$ (its shape becomes independent of $t$).\\

\noindent
Interestingly, Eq.~\eqref{eq:main_cav} implies that we can use $h_t^{\alpha}(\tau)$ to convert Eq.~\eqref{eq:cav} into the following Riemann integral:
\begin{equation}
\label{eq:option2}
\int_0^{\frac{t^{\alpha}}{\Gamma(\alpha+1)}} f\circ h_{\alpha}^t (\tau) d\tau.  
\end{equation}
We can, therefore, evaluate Eq.~\eqref{eq:cav} in one of two ways: we can either evaluate Eq.~\eqref{eq:rs} or we can evaluate Eq.~\eqref{eq:option2}.

\section{Examples}
In this section we will illustrate the usefulness of the non-static non-rectangular sum-based geometric interpretation we presented in the previous section at the 
hand of two examples.\\

\begin{figure}[htb]
\centering
\includegraphics[width=0.75\textwidth]{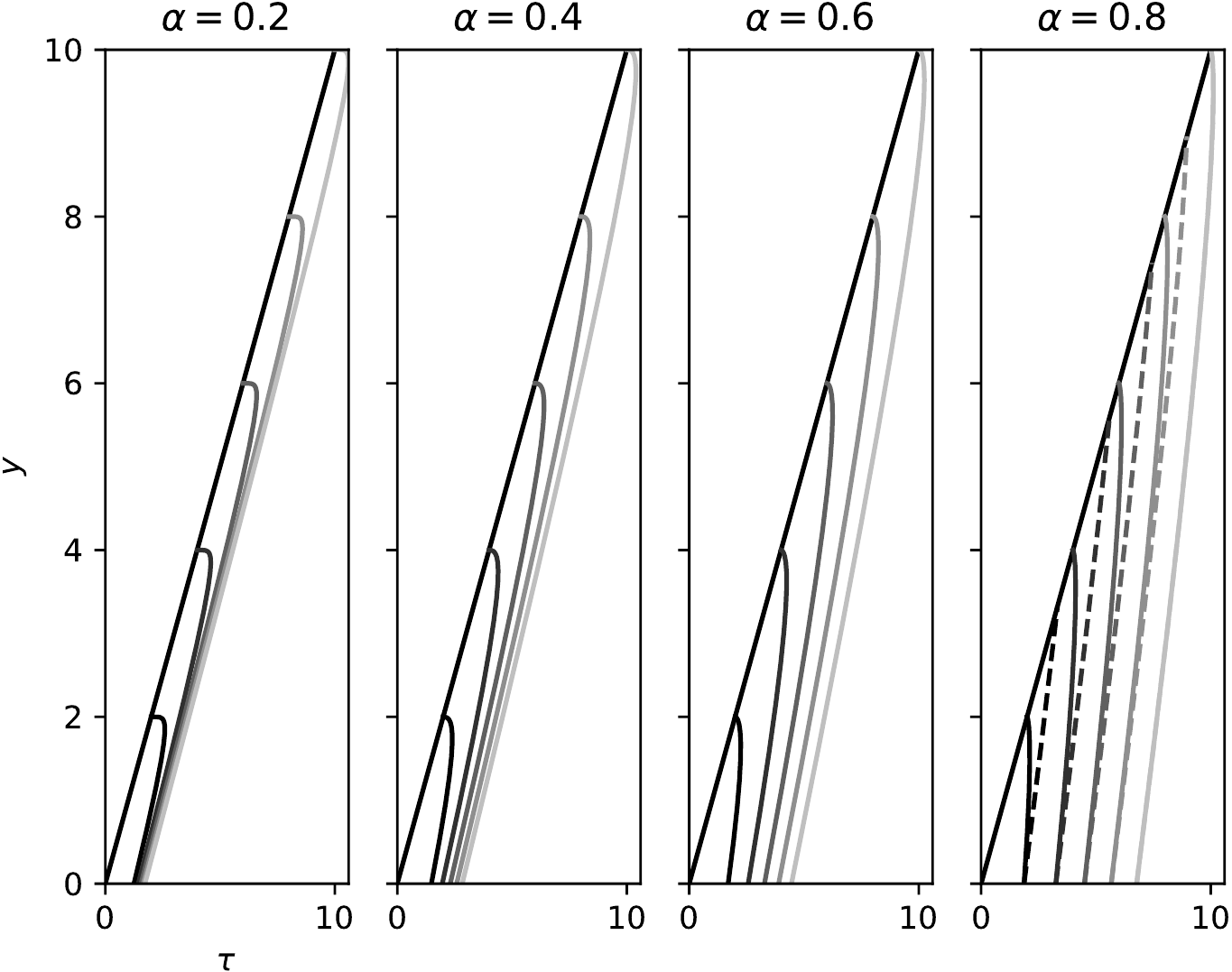}
\caption{The non-static non-rectangular sum-based geometric interpretation of Eq.~\eqref{eq:ex1_raw}.}
\label{fig:geo1}
\end{figure}

\noindent
Let us first consider the following fractional integral:
\begin{equation}
\label{eq:ex1_raw}
\frac{1}{\Gamma(\alpha)}\int_0^t f(\tau) (t-\tau)^{\alpha-1}~d\tau = \frac{1}{\Gamma(\alpha)}\int_0^t \tau(t-\tau)^{\alpha-1}~d\tau. 
\end{equation}
The function $f(\tau) =\tau$ is depicted in Fig.~\ref{fig:geo1}. As we mentioned in the previous section, we can assign a geometric interpretation to the above fractional integral by converting it into a Cavalieri integral using two steps. First, we convert 
the above integral into a Riemann-Stieltjes integral:
\begin{equation}
\frac{1}{\Gamma(\alpha)}\int_0^t f(\tau) (t-\tau)^{\alpha-1}~d\tau = \int_0^t f(\tau)~dg_t^{\alpha}(\tau). 
\end{equation}
We then convert the resulting Riemann-Stieltjes integral into a Cavalieri integral:
\begin{equation}
\label{eq:ex1}
\frac{1}{\Gamma(\alpha)}\int_0^t f(\tau) (t-\tau)^{\alpha-1}~d\tau = \int_0^t f(\tau)~dg_t^{\alpha}(\tau)=\int_{a_t^{\alpha}(y)}^{b_t^{\alpha}(y)} f(\tau)~d\tau,
\end{equation}
with $a_t^{\alpha}(y) = y - \frac{t^{\alpha}-(t-y)^{\alpha}}{\Gamma(\alpha+1)}.$

\noindent
Since Eq.~\eqref{eq:ex1}, is a Cavalieri integral it has a known geometric interpretation (see Fig.~\ref{fig:caval2}). It can be interpreted as the area obtained by summing together the areas of an infinite 
number of infinitesimals whose shape is determined by $\alpha$ and $t$, i.e. these infinitesimals have to be translations of $b_t^{\alpha}(y)$. The geometric interpretation of Eq.~\eqref{eq:ex1_raw} is, therefore, non-static as it depends on the value of $\alpha$ and $t$. 
There is one thing, however, which remains unclear: in what way do $\alpha$ and $t$ affect the shape of the aforementioned infinitesimals, i.e. how does the value of $\alpha$ and $t$ affect the shape of $b_t^{\alpha}(y)$. 
We will answer this question by working through a detailed example.\\


\noindent
If we set the value of $\alpha$ to $\sfrac{4}{5}$ and the value of $t$ to 10 in Eq.~\eqref{eq:ex1} we obtain:
\begin{equation}
\label{eq:ex1_specific}
A_{10}^{\sfrac{4}{5}}=\int_{a_{10}^{\sfrac{4}{5}}(y)}^{b_{10}^{\sfrac{4}{5}}(y)} f(\tau)~d\tau 
\end{equation}
Let $R_{10}^{\sfrac{4}{5}}$ denote the region which is bounded by the $\tau$-axis, the function $f(\tau)$ and the function $b_{10}^{\sfrac{4}{5}}(y)$. The region $R_{10}^{\sfrac{4}{5}}$ is depicted in the right-most panel of Fig.~\ref{fig:geo1}. It is the largest of all the 
regions that are depicted in the right-most panel of Fig.~\ref{fig:geo1}. Moreover, let $A_{10}^{\sfrac{4}{5}}$ denote the area of $R_{10}^{\sfrac{4}{5}}$. 
We can approximate the area of $R_{10}^{\sfrac{4}{5}}$ by summing together the area of $n$ equal width non-rectangular integration strips that inscribe $R_{10}^{\sfrac{4}{5}}$, i.e. we can 
construct a Cavalieri sum (see Eq.~\eqref{eq:cav_sum}). Note, the sides of the non-rectangular integration strips we use in this approximation have to be translations of $b_{10}^{\sfrac{4}{5}}(y)$.
If we take the limit of the aforementioned sum as $n\rightarrow \infty$ we approach the Cavalieri integral in Eq.~\eqref{eq:ex1_specific}, i.e. $A_{10}^{\sfrac{4}{5}}$.  
Eq.~\eqref{eq:ex1_specific}, therefore, implies that the region $R_{10}^{\sfrac{4}{5}}$ is made up 
of an infinite number of non-rectangular infinitesimals; all of them being translations of $b_{10}^{\sfrac{4}{5}}(y)$. Note, the right-most curve depicted in the right-most panel of Fig.~\ref{fig:geo1} is $b_{10}^{\sfrac{4}{5}}(y)$.  Moreover, the special case in which we approximate $A_{10}^{\sfrac{4}{5}}$ using five equal width non-rectangular integration strips is depicted in the right-most panel of Fig.~\ref{fig:geo1}. 
The sides of the aforementioned five integration strips are depicted in Fig.~\ref{fig:geo1} using dashed lines. Let, $\mathcal{B}_{\mathsf{T}}^{\sfrac{4}{5}}=\{b_{10}^{\sfrac{4}{5}}(y)-c|c\in(\mathsf{T}-2)\}$ with $\mathsf{T}=\{2,4,6,8,10\}$.\\ 


\noindent
Varying both $\alpha$ and $t$ results in Eq.~\eqref{eq:ex1}. In contrast with Eq.~\eqref{eq:ex1_specific}, Eq.~\eqref{eq:ex1} can be associated with more than one region. 
Let $R_{t}^{\alpha}$ denote the region which is bounded by the $\tau$-axis, the function $f(\tau)$ and the function $b_{t}^{\alpha}(y)$. Note that $b_{t}^{\alpha}(y)$ forms the 
right-most edge of $R_{t}^{\alpha}$. Moreover, let $A_{t}^{\alpha}$ denote the area of $R_{t}^{\alpha}$ and let $\mathsf{R}_{\mathsf{T}}^{\mathsf{A}}=\{R_t^{\alpha}|t\in\mathsf{T},~\alpha\in\mathsf{A}\}$ and 
$\mathsf{B}_{\mathsf{T}}^{\mathsf{A}}=\{b_t^{\alpha}(y)|t\in\mathsf{T},~\alpha\in\mathsf{A}\}$. 
The regions in the set $\mathsf{R}_{\mathsf{T}}^{\mathsf{A}}$ are depicted in Fig.~\ref{fig:geo1}. The regions that are depicted in each of the panels of Fig.~\ref{fig:geo1} were generated by varying $t$, whilst keeping the value of $\alpha$ fixed. 
The title of each panel indicates which $\alpha$-value was used to generate the regions depicted in each panel. The only difference between the regions located in each 
panel is their right edges. The right edges of the regions in each panel are depicted using varying shades of gray. The lighter the shade with which an edge is depicted, 
the larger the $t$-value is that is associated with it.\\

\noindent
As was the case for Eq.~\eqref{eq:ex1_specific} and $R_{10}^{\sfrac{4}{5}}$, Eq.~\eqref{eq:ex1} implies that we can assume that each of the regions depicted in Fig.~\ref{fig:geo1} are made up of an infinite number of infinitesimals. The infinitesimals that make up a specific region are translations of that region's right edge. 
The following observations can be made by inspecting Fig.~\ref{fig:geo1}: 
\begin{description}
 \item[Effect of $\alpha$] the gray-tinted solid curves depicted in the different panels of Fig.~\ref{fig:geo1} are not equal to one another, i.e. $\mathsf{B}_{\mathsf{T}}^{\alpha_1}\neq\mathsf{B}_{\mathsf{T}}^{\alpha_2}$. 
 \item[Effect of $t$] the gray-tinted solid curves depicted in each of the panels of Fig.~\ref{fig:geo1} are not translations of one another. Moreover, the dashed gray curves and the solid gray-tinted curves in the right-most panel are not equal to one another, i.e. $\mathsf{B}_{\mathsf{T}}^{\sfrac{4}{5}}\neq\mathcal{B}_{\mathsf{T}}^{\sfrac{4}{5}}$.
\end{description}
We can therefore draw the following conclusion: the value of both $\alpha$ and $t$ alters the shape of the infinitesimals ($b_{t}^{\alpha}(y)$ depends on $\alpha$ and $t$) that make up the regions in Fig.~\ref{fig:geo1}. Interestingly, the effect of $\alpha$ on the shape of the infinitesimals that make up the regions depicted in Fig.~\ref{fig:geo1} is larger than the effect that $t$ has on their shape. 
The curves associated with the infinitesimals obtained by varying $t$, whilst keeping the value of $\alpha$ fixed are very similar in nature and there exists a high degree of correlation between them (if their 
sizes are ignored). The curves associated with the infinitesimals obtained by varying $\alpha$, whilst keeping the value of $t$ fixed are quite different from one another (the aforementioned curves are not highly correlated with 
one another).\\  


\begin{table}[h!]
\centering
\caption{Curves obtained by evaluating Eq.~\eqref{eq:ex1} and Eq.~\eqref{eq:ex2} for all $\alpha\in \mathsf{A}$ using either Eq.~\eqref{eq:rs} or Eq.~\eqref{eq:option2}.}
\label{tab:eval}
\begin{tabular}{||c|| c| c||} 
 \hline
 $\alpha$ & $\frac{1}{\Gamma(\alpha)}\int_0^t \tau(t-\tau)^{\alpha-1}~d\tau$ & $\frac{1}{\Gamma(\alpha)}\int_0^t \sqrt{\tau}(t-\tau)^{\alpha-1}~d\tau$ \\[0.5ex]
 \hline\hline
 \rule{0pt}{2.5ex}
 $0$ & $t$ & $\sqrt{t}$  \\  
 \rule{0pt}{2.5ex}
 $0.2$ & $\frac{25}{6\Gamma(\frac{1}{5})} t^{\sfrac{6}{5}}$  & $\frac{\sqrt{\pi}}{2\Gamma(\frac{17}{10})} t^{\sfrac{7}{10}}$  \\ 
 \rule{0pt}{2.5ex}
 $0.4$ & $\frac{25}{14\Gamma(\frac{2}{5})} t^{\sfrac{7}{5}}$  & $\frac{\sqrt{\pi}}{2\Gamma(\frac{19}{10})} t^{\sfrac{9}{10}}$  \\
 \rule{0pt}{2.5ex}
 $0.6$ & $\frac{25}{24\Gamma(\frac{3}{5})} t^{\sfrac{8}{5}}$  & $\frac{\sqrt{\pi}}{2\Gamma(\frac{21}{10})} t^{\sfrac{11}{10}}$  \\ 
 \rule{0pt}{2.5ex}
 $0.8$ & $\frac{25}{36\Gamma(\frac{4}{5})} t^{\sfrac{9}{5}}$  & $\frac{\sqrt{\pi}}{2\Gamma(\frac{23}{10})} t^{\sfrac{13}{10}}$  \\ 
 \rule{0pt}{2.5ex}
 $1$ & $\frac{1}{2}t^2$  & $\frac{2}{3}t^{\sfrac{3}{2}}$  \\ [1ex]
 \hline
 \end{tabular}
\end{table}

\begin{figure}[htb]
\centering
\includegraphics[width=0.75\textwidth]{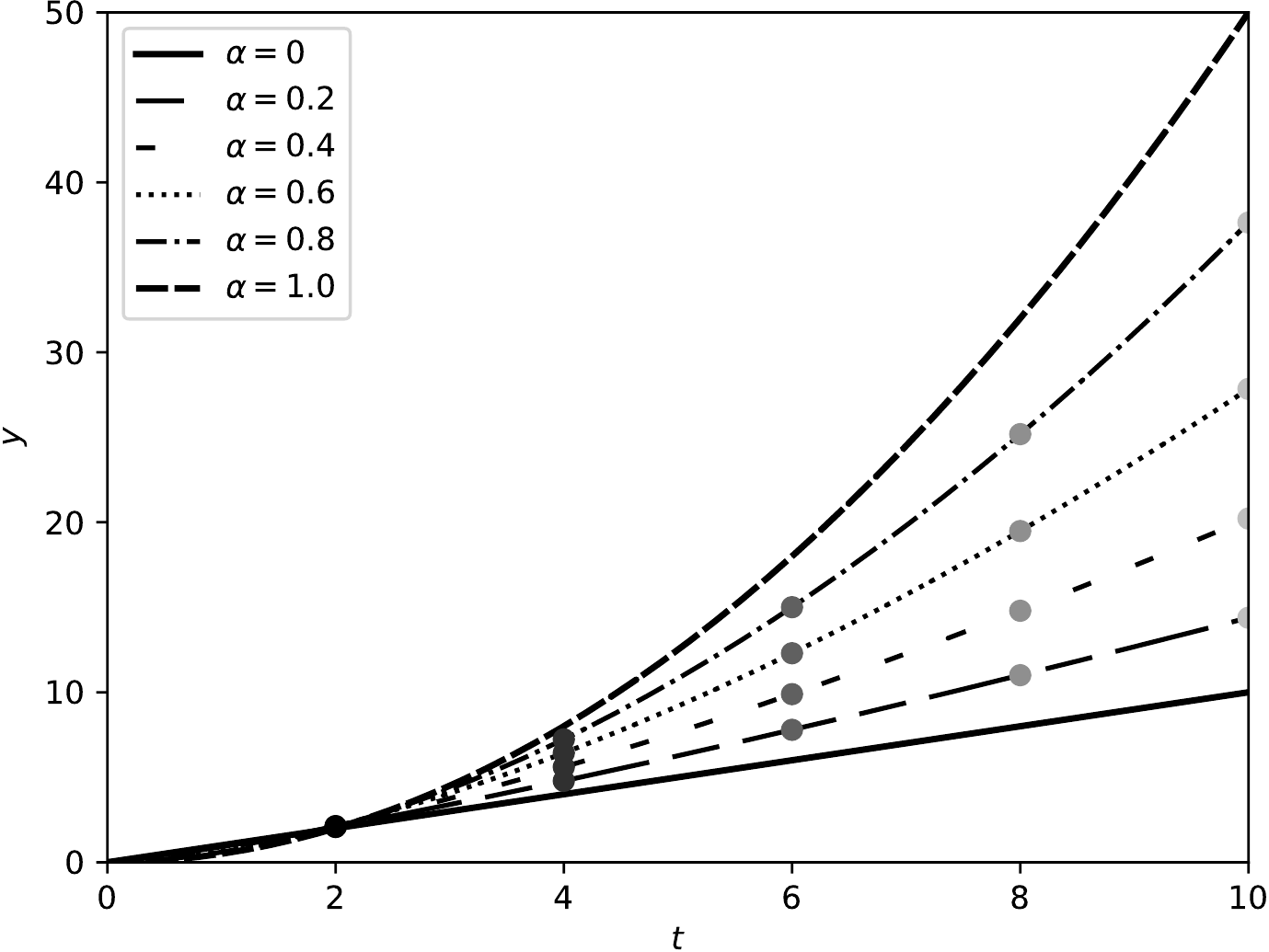}
\caption{Depicts the curves which are obtained by evaluating Eq.~\eqref{eq:ex1} for all $\alpha\in \mathsf{A}$ using either Eq.~\eqref{eq:rs} or Eq.~\eqref{eq:option2}.}
\label{fig:eval1}
\end{figure}

\noindent
If we evaluate the integral in Eq.~\eqref{eq:ex1} for all $\alpha\in \mathsf{A}$ (using either Eq.~\eqref{eq:rs} or Eq.~\eqref{eq:option2}) we obtain the curves in the 
left-most column of Table~\ref{tab:eval}. These curves are depicted in Fig.~\ref{fig:eval1}.
The area of the regions depicted in Fig.~\ref{fig:geo1} are also plotted in Fig.~\ref{fig:eval1} using circular markers. As expected, the circular markers in Fig.~\ref{fig:eval1} that correspond to the areas of the regions in Fig.~\ref{fig:eval1}  
associated with a particular choice of $\alpha$ (i.e. regions depicted in a particular panel of Fig.~\ref{fig:geo1}) fall on the curve in Fig.~\ref{fig:eval1}; associated with the same choice of $\alpha$. Moreover, note that we have adopted the 
same coloring scheme in Fig.~\ref{fig:geo1} and Fig.~\ref{fig:eval1}. Using the same coloring scheme in both figures makes it clear which of the area values in Fig.~\ref{fig:eval1} can be associated with which regions in Fig.~\ref{fig:geo1}.
To summarize, Fig.~\ref{fig:geo1}, therefore, depicts a non-static non-rectangular sum-based geometric interpretation of Eq.~\eqref{eq:ex1}. Moreover, Fig.~\ref{fig:eval1} corroborates this 
geometric interpretation.\\

\begin{figure}[htb]
\centering
\includegraphics[width=0.75\textwidth]{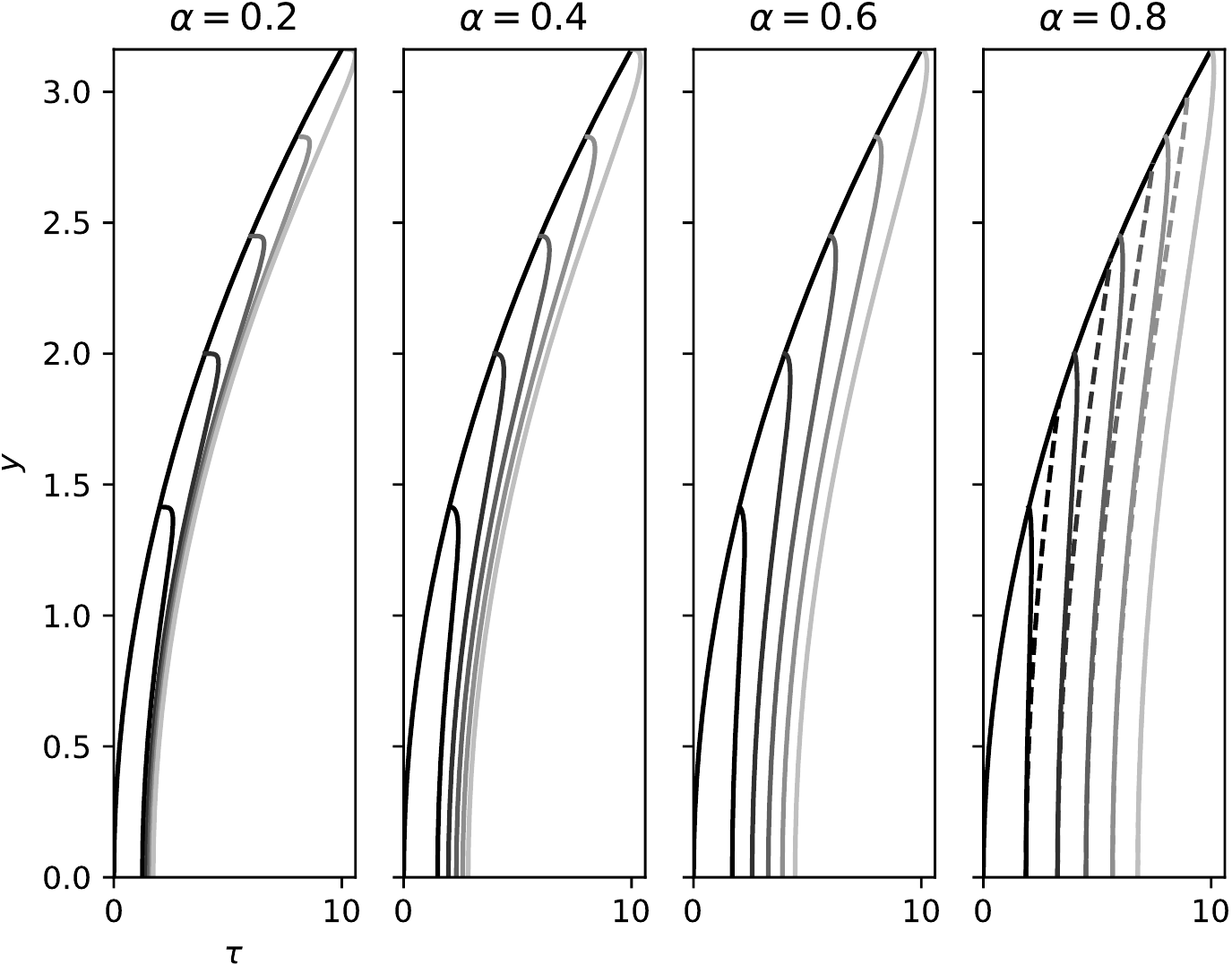}
\caption{The non-static non-rectangular sum-based geometric interpretation of Eq.~\eqref{eq:ex2}.}
\label{fig:geo2}
\end{figure}

\begin{figure}[htb]
\centering
\includegraphics[width=0.75\textwidth]{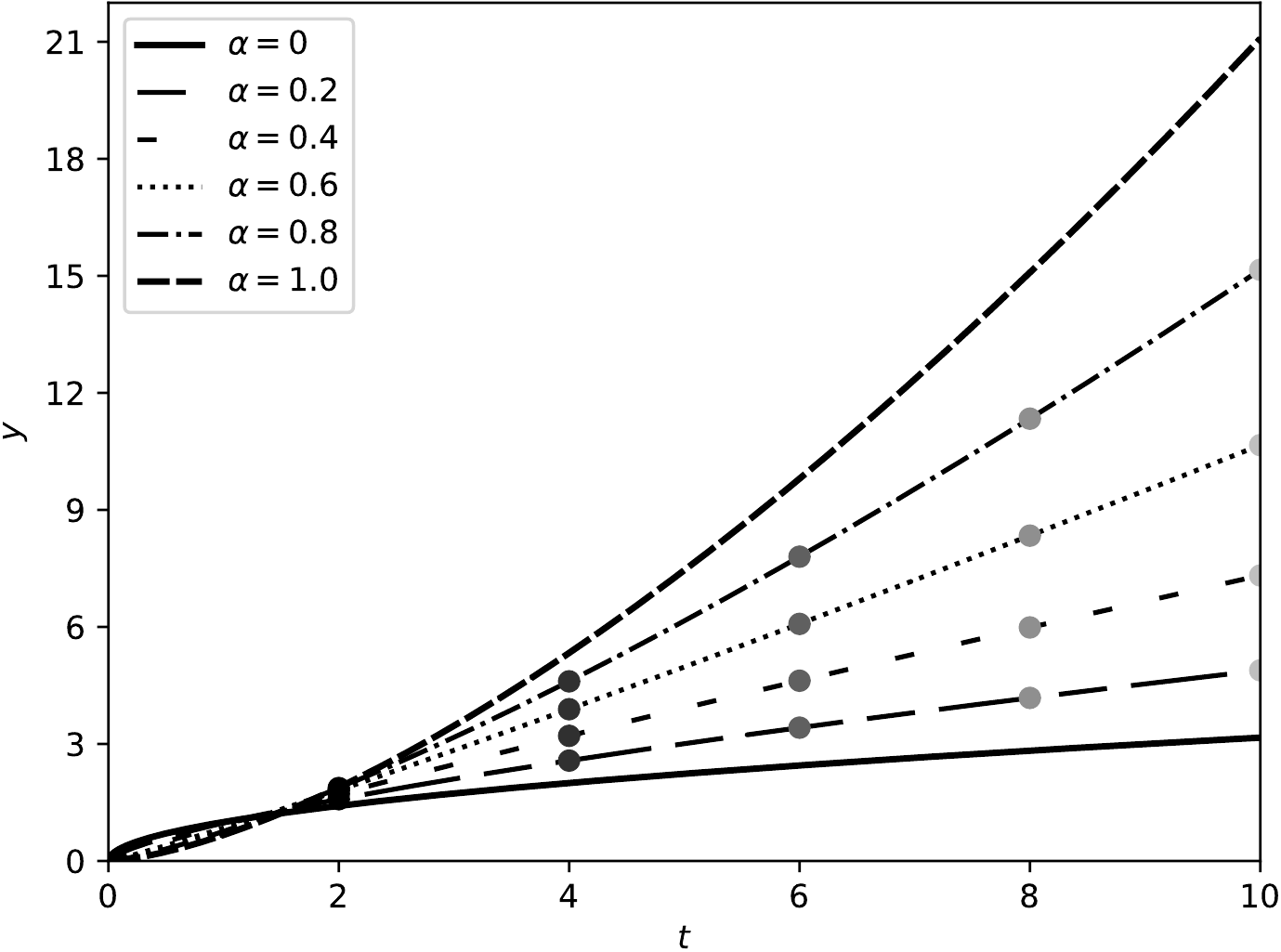}
\caption{Depicts the curves which are obtained by evaluating Eq.~\eqref{eq:ex2} for all $\alpha\in \mathsf{A}$ using either Eq.~\eqref{eq:rs} or Eq.~\eqref{eq:option2}.}
\label{fig:eval2}
\end{figure}

\noindent
Let us now consider a different fractional integral:
\begin{equation}
\label{eq:ex2}
\frac{1}{\Gamma(\alpha)}\int_0^t f(\tau) (t-\tau)^{\alpha-1}~d\tau = \frac{1}{\Gamma(\alpha)}\int_0^t \sqrt{\tau}(t-\tau)^{\alpha-1}~d\tau. 
\end{equation}
The non-static non-rectangular sum-based geometric interpretation associated with Eq.~\eqref{eq:ex1} is depicted in Fig.~\ref{fig:geo2}. 
If we evaluate the integral in Eq.~\eqref{eq:ex1} for all $\alpha\in \mathsf{A}$ we obtain the curves in the right-most column of Table~\ref{tab:eval}. 
These curves are depicted in Fig~\ref{fig:eval2}.
The validity of the non-static non-rectangular sum-based geometric interpretation depicted in Fig.~\ref{fig:geo2} is corroborated by Fig.~\ref{fig:eval2}. 

\vfill\eject

\end{document}